\documentclass[11pt]{article}

\usepackage{latexsym}
\usepackage{amsfonts}
\usepackage{amsmath}
\usepackage{amsthm}
\usepackage{mathrsfs}  
\usepackage{dsfont}    
\usepackage{bbold}     
\usepackage[english]{babel}
\usepackage{caption}
\usepackage{epsfig}
\usepackage{float}
\usepackage{subfigure}

\DeclareMathOperator{\Val}{Val}

\newtheorem{theorem}{Theorem}
\newtheorem*{theorem*}{Theorem}

\newtheorem{lemma}[theorem]{Lemma}

\newtheorem{hyp}{Assumption}

\newcommand{\zerarcounters}{\setcounter{equation}{0}\setcounter{theorem}{0}}

\newcommand{\ZZZ}{\mathds{Z}}

\newcommand{\NNN}{\mathds{N}}

\newcommand{\RRR}{\mathds{R}}
\newcommand{\TTT}{\mathds{T}}

\newcommand{\MM}{{\mathcal M}}

\newcommand{\TT}{{\mathcal T}}

\newcommand{\gotn}{{\mathfrak n}}

\newcommand{\gotp}{{\mathfrak p}}

\newcommand{\gotB}{{\mathfrak B}}

\newcommand{\gotN}{{\mathfrak N}}

\newcommand{\gotR}{{\mathfrak R}}

\newcommand{\gotT}{{\mathfrak T}}

\newcommand{\Fullbox}{{\rule{2.0mm}{2.0mm}}}

\newcommand{\EP}{\hfill\Fullbox\vspace{0.2cm}}
\newcommand{\prova}{\noindent{\it Proof. }}
\newcommand{\io}{\infty}
\newcommand{\e}{\varepsilon}
\newcommand{\al}{\alpha}

\newcommand{\g}{\gamma}
\newcommand{\om}{\omega}

\newcommand{\laa}{\langle}
\newcommand{\raa}{\rangle}

\newcommand{\oo}{\boldsymbol{\omega}}
\newcommand{\nn}{\boldsymbol{\nu}}
\newcommand{\pps}{\boldsymbol{\psi}}
\newcommand{\vzero}{\boldsymbol{0}}
\newcommand{\der}{{\rm d}}
\newcommand{\ii}{{\rm i}}

\def\ins#1#2#3{\vbox to0pt{\kern-#2 \hbox{\kern#1 #3}\vss}\nointerlineskip}

\begin{document}

\title{\bf Construction of quasi-periodic response solutions\\
in forced strongly dissipative systems}
\author
{\bf Guido Gentile
\vspace{2mm}
\\ \small 
Dipartimento di Matematica, Universit\`a di Roma Tre, Roma,
I-00146, Italy.
\\ \small 
E-mail: gentile@mat.uniroma3.it}

\date{}

\maketitle

\begin{abstract}
We consider a class of ordinary differential equations describing
one-dimensional quasi-periodically forced systems in the presence
of large damping. We give a fully constructive proof of the existence
of response solutions, that is quasi-periodic solutions which have
the same frequency vector as the forcing. This requires dealing with
a degenerate implicit function equation: we prove that the latter
has a unique solution, which can be explicitly determined.
As a by-product we obtain an explicit estimate of the
minimal size of the damping coefficient.
\end{abstract}




\zerarcounters
\section{Introduction}
\label{sec:1}

In this paper we continue the analysis started in \cite{GBD1,GBD2,G2}
on the existence and properties of quasi-periodic motions
in one-dimensional strongly dissipative forced systems.

We consider one-dimensional systems with a quasi-periodic forcing
term in the presence of strong damping, described
by ordinary differential equations of the form
\begin{equation}
\e \ddot x + \dot x + \e g(x) = \e f(\oo t) ,
\label{eq:1.1} \end{equation}
where $g\!:\RRR\to\RRR$ and $f\!:\TTT^{d}\to\RRR$ are 
real analytic functions and $\TTT=\RRR/2\pi\ZZZ$.
We call $g(x)$ the \textit{mechanical force},
$f(\oo t)$ the \textit{forcing term},
$\oo\in\RRR^{d}$ the \textit{frequency vector} of the forcing,
and $\g=1/\e>0$ the \textit{damping coefficient}.

The function $f$ is quasi-periodic in $t$, i.e.
\begin{equation}
f(\pps) = \sum_{\nn\in \ZZZ^{d}}
{\rm e}^{\ii\nn \cdot \pps} f_{\nn} , \qquad \pps \in \TTT^{d} ,
\label{eq:1.2} \end{equation}
with average $\laa f \raa = f_{\vzero}$, and
$\cdot$ denoting the scalar product in $\RRR^{d}$.
By the analyticity assumption on $f$ and $g$, one has
$|f_{\nn}| \le \Phi {\rm e}^{-\xi|\nn|}$ and
$|\der^{s}g(c_{0})/\der x^{s}| \le s!\Gamma^{s}$
for suitable positive constants $\Phi$, $\xi$, and $\Gamma$.

A Diophantine condition is assumed on $\oo$. Define
the \textit{Bryuno function}
\begin{equation}
\gotB(\oo) = \sum_{n=0}^{\infty} \frac{1}{2^{n}} \log
\frac{1}{\al_{n}(\oo)} , \qquad \al_{n}(\oo) = \inf
\{ |\oo\cdot\nn| : \nn \in \ZZZ^{d}
\hbox{ such that } 0<|\nn|\le 2^{n} \} .
\label{eq:1.3} \end{equation}
%
\begin{hyp}\label{hyp:1} 
The frequency vector $\oo$ satisfies
the Bryuno condition $\gotB(\oo)<\infty$.
\end{hyp}

The following assumption will be made on the functions $g$ and $f$ (for
given force $g(x)$ this can read as a condition on the forcing term).

\begin{hyp}\label{hyp:2} 
There exists $c_{0}\in\RRR$ such that $x=c_{0}$ is a zero
of odd order $\gotn$ of the equation
\begin{equation}
g(x) - f_{\vzero} = 0 ,
\label{eq:1.4} \end{equation}
that is $g_{0}=\gotn!^{-1}
\der^{\gotn}g/\der x^{\gotn}(c_{0})\neq0$ and, if $\gotn>1$,
$\der^{k} g/\der x^{k}(c_{0})=0$ for $k=1,\ldots,\gotn -1$.
\end{hyp}

In \cite{G2} we proved the following result about the existence
of quasi-periodic solutions with the same frequency vector $\oo$
as the forcing (response solutions).

\begin{theorem} \label{thm:1.1}
Under Assumptions \ref{hyp:1} and \ref{hyp:2}, for $\e$
small enough there exists at least one quasi-periodic solution
$x_{0}(t)=c_{0}+u(\oo t,\e)$, with frequency vector $\oo$,
reducing to $c_{0}$ as $\e$ tends to $0$.
Such a solution is analytic in $t$.
\end{theorem}

Note that the condition $\e>0$ could be eliminated:
indeed, the proof in \cite{G2} works for all $\e$ small enough,
and the request $\e>0$ only aims to interpret $\g=1/\e$
as the damping coefficient. Analyticity on $t$ is not explicitly
stated in \cite{G2}, but follows immediately from the proof therein.

We also proved in \cite{G2} that the condition that $c_{0}$ be a zero of
odd order of (\ref{eq:1.4}) is a necessary and sufficient condition for
a quasi-periodic solution around $c_{0}$ to exists.

However, as pointed out in \cite{G2}, except for the non-degenerate
case $\gotn=1$ (where the implicit function theorem applies),
in general the proof ultimately relies on continuity arguments,
which do not provide a quantitative constructive estimate
on the maximal size of the perturbation parameter $\e$. Moreover,
the quasi-periodic solution was constructed in terms of two parameters,
that is $\e$ and $c$: the latter is defined as the constant part of
the quasi-periodic solution itself and is fixed in terms of $\e$
so as to solve a certain implicit function equation (the so-called
bifurcation equation). In particular, a quasi-periodic solution
was showed to exist for any solution $c(\e)$ to the bifurcation equation,
but the problem of studying how many such solutions exist and
how do they depend on $\e$ was not investigated.

In this paper we show that for all odd $\gotn$ the bifurcation equation
admits one and only one solution. Furthermore, we explicitly construct
the quasi-periodic solution in terms of the only parameter $\e$,
and we also give a quantitative estimate on the maximal size of $\e$.
So, with respect to \cite{G2}, the proof of existence
of the quasi-periodic solution is fully constructive.

We can summarise our results in the following statement.

\begin{theorem} \label{thm:1.2}
Under Assumptions \ref{hyp:1} to \ref{hyp:2} on the
ordinary differential equation (\ref{eq:1.1}), there exists
an (explicitly computable) constant $\e_{0}>0$ such that
for all $|\e|<\e_{0}$ there exists a quasi-periodic solution
$x_{0}(t)=c_{0}+u(\oo t,\e)$, with frequency vector $\oo$,
such that $u=O(\e)$ as $\e\to0$. Such a solution is analytic
in $t$ and depends $C^{\io}$-smoothly on $\e$.
\end{theorem}

The result solves a problem left as open in \cite{G2}.
The problem remains whether other quasi-periodic solutions exist.

Numerical algorithms to construct response solutions of
quasi-periodically forced dissipative systems are provided in
\cite{CH,KP1,KP2}, based either on a generalised harmonic balance
method or on a fixed point method for a suitable Poincar\'e map.
Also the method described in the present paper is well suited
for numerical implementations, and allows a completely rigorous
control of the approximation error for the solution.
Moreover it applies also in degenerate cases where one cannot
apply directly the implicit function theorem (of course
it requires for the damping coefficient to be large enough).

\zerarcounters
\section{The bifurcation equation}
\label{sec:2}

The existence of quasi-periodic solutions $x_{0}(t)$ has been proved in
\cite{G2}. The proof proceeds as follows. First of all, write
\begin{equation}
x_{0}(t)=c_{0}+u(\oo t,\e)=c+X(\oo t;\e,c) , \qquad
X(\pps;\e,c)=\sum_{\nn\in\ZZZ^{d}_{*}} {\rm e}^{\ii\nn\cdot\pps}
X_{\nn}(\e,c) ,
\label{eq:2.1} \end{equation}
with $\ZZZ^{d}_{*}=\ZZZ^{d}\setminus\{\vzero\}$ (so that $c$ is the
average of $x_{0}$ on the $d$-dimensional torus).
The Fourier coefficients $X_{\nn}=X_{\nn}(\e,c)$ of the
function $X(\pps;\e,c)$ are obtained by solving
the \textit{range equation}
\begin{equation}
\ii\oo\cdot\nn \left( 1 + \ii\e\oo\cdot\nn \right) X_{\nn} +
\e \left[ g(c+X(\cdot;\e,c) \right]_{\nn} = \e f_{\nn} ,
\qquad \nn \neq \vzero ,
\label{eq:2.2} \end{equation}
where
\begin{equation}
\left[ g(c+X(\cdot;\e,c) \right]_{\nn} =
\sum_{p=0}^{\infty} \frac{1}{p!} \frac{{\rm d}^{p}}{{\rm d}x^{p}}
g(c) \sum_{\substack{\nn_{1},\ldots,\nn_{p}\in\ZZZ^{d}_{*} \\
\nn_{1}+\ldots+\nn_{p}=\nn}} X_{\nn_{1}}(\e,c) \ldots
X_{\nn_{p}}(\e,c) .
\label{eq:2.3} \end{equation}
The analysis in \cite{G2} shows that for all $c$ close enough
to $c_{0}$ and all $\e$ small enough there exists a function
$X(\pps;\e,c)$ which solves (\ref{eq:2.2}). Moreover
the map $(\e,c) \mapsto X(\cdot;\e,c)$
is $C^{\infty}$ in a neighbourhood of $(0,c_{0})$, and
\begin{equation}
X(\pps;\e,c) = \e X^{(1)}(\pps) + o(\e) , \qquad
\dot X^{(1)}(\oo t) = f(\oo t) .
\label{eq:2.4} \end{equation}
Then, for any solution $c=c(\e)$ to the \textit{bifurcation equation}
\begin{equation}
F(\e,c) := [g(c+X(\cdot;\e,c)]_{\vzero} - f_{\vzero} = 0 ,
\label{eq:2.5} \end{equation}
there exists a quasi-periodic solution $x_{0}(t)=c(\e)+X(\oo t;\e,c(\e))$.

In principle, there could be several solutions to (\ref{eq:2.5})
when $\gotn>1$. On the contrary, we shall prove that
the solution to (\ref{eq:2.5}) is unique.
Unfortunately, this does not implies that the
response solution $c_{0}+u(\oo t,\e)$ of Theorem \ref{thm:1.1}
is the only quasi-periodic solution reducing to $c_{0}$ as $\e\to0$,
because no result ensures that the function $X(\pps;\e,c)$
which solves (\ref{eq:2.2}) at fixed $c$ and $\e$ is unique.

We shall write $c=c(\e)=c_{0}+\zeta$, with $\zeta=\zeta(\e)$
such that $\zeta(0)=0$, so that (\ref{eq:2.3}) becomes
\begin{equation}
\left[ g(c+X(\cdot;\e,c) \right]_{\nn} =
\sum_{p=\gotn}^{\infty} \frac{1}{p!} \frac{{\rm d}^{p}}{{\rm d}x^{p}}
g(c_{0}) \sum_{k=0}^{p} 
\left( \begin{matrix} p \\ k \end{matrix} \right) \zeta^{p-k}
\sum_{\substack{\nn_{1},\ldots,\nn_{k}\in\ZZZ^{d}_{*} \\
\nn_{1}+\ldots+\nn_{k}=\nn}} X_{\nn_{1}}(\e,c) \ldots
X_{\nn_{k}}(\e,c) .
\label{eq:2.6} \end{equation}
%

\begin{lemma} \label{lem:2.1}
Let $c=c(\e)$ a solution of the bifurcation equation (\ref{eq:2.5}).
Then $c=c_{0}+O(\e)$, and $c=c_{0}+o(\e)$
requires $[(X^{(1)}(\cdot))^{\gotn}]_{\vzero}=0$.
\end{lemma}

\prova Assume that $\zeta/\e\to\io$ as $\e\to0$.
Then, by using (\ref{eq:2.4}), one has
\begin{equation}
[g(c+X(\cdot;\e,c)]_{\vzero} - f_{\vzero} =
g_{0} \left[ \left( \zeta + X(\cdot,\e,c)
\right)^{\gotn} \right]_{\vzero} +
O(\zeta^{\gotn+1}) =
g_{0} \zeta^{\gotn} + O(\zeta^{\gotn+1}) = 0 ,
\nonumber \end{equation}
which leads to a contradiction. Hence $\zeta=O(\e)$.
On the other hand if $[(X^{(1)}(\cdot))^{\gotn}]_{\vzero}\neq0$
and $\zeta=o(\e)$ one has
\begin{equation}
[g(c+X(\cdot;\e,c)]_{\vzero} - f_{\vzero} =
g_{0} \left[ \left( \zeta +
X(\cdot;\e,c) \right)^{\gotn} \right]_{\vzero} +
O(\e^{\gotn+1}) =
g_{0} \e^{\gotn} [(X^{(1)}(\cdot))^{\gotn}]_{\vzero} +
O(\e^{\gotn+1}) = 0 ,
\nonumber \end{equation}
which once more leads to a contradiction.\EP

Note that a sufficient condition for
$[(X^{(1)}(\cdot))^{\gotn}]_{\vzero}$
to vanish is that $f$ in (\ref{eq:1.2}) is even.

\begin{figure}[ht]
\vskip.3truecm
\centering
\ins{031pt}{2pt}{$j$}
\ins{138pt}{-76pt}{$k$}
\ins{176pt}{2pt}{$j$}
\ins{282pt}{-76pt}{$k$}
\ins{321pt}{2pt}{$j$}
\ins{428pt}{-76pt}{$k$}
\ins{030pt}{-13pt}{$\gotn$}
\ins{175pt}{-13pt}{$\gotn$}
\ins{320pt}{-13pt}{$\gotn$}
\ins{094pt}{-78pt}{$\gotn$}
\ins{220pt}{-76pt}{$\gotn\!-\!1$}
\ins{364pt}{-76pt}{$\gotn\!-\!1$}
\ins{177pt}{-56pt}{$1$}
\ins{321pt}{-56pt}{$1$}
\ins{412pt}{-76pt}{$\gotp$}
\includegraphics{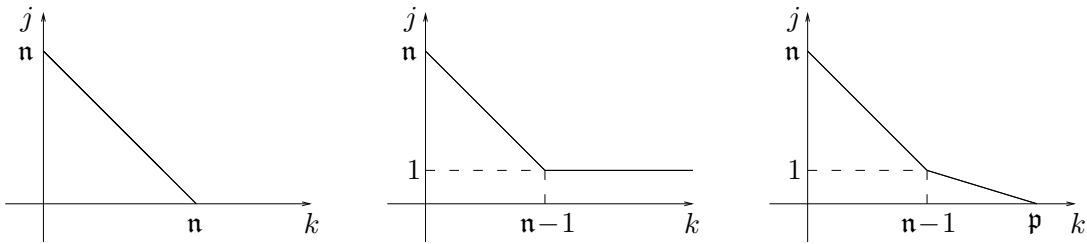}
\caption{Newton polygon (corresponding to cases 1 to 3, respectively)}
\label{fig:1}
\end{figure}

The function $F(\e,c)$, defined in (\ref{eq:2.5}),
is $C^{\io}$ in both $\e$ and $c$; see \cite{G2}.
In order to identify the leading orders to $F(\e,c)$,
we consider the carrier
\begin{equation}
\Delta(F) = \left\{ (k,j) \in \ZZZ_{+}\times\ZZZ_{+} :
F_{k,j} \neq 0 \right\} , \qquad
F_{k,j} = \frac{{\rm d}^{k}}{{\rm d}\e^{k}}
\frac{{\rm d}^{j}}{{\rm d}c^{j}}F(0,c_{0}) ,
\label{eq:2.7} \end{equation}
and draw the Newton polygon in the $(k,j)$ plane; see \cite{BK}
(see also \cite{CG}). If we denote by $\{(k_{1},j_{1}),(k_{2},j_{2})\}$
the segment with endpoints $(k_{1},j_{1})$ and $(k_{2},j_{2})$
in the $(k,j)$ plane, by construction for $\gotn>1$
the Newton polygon consists of (cf. figure \ref{fig:1})
\begin{enumerate}
\item either the only segment $\{(0,\gotn),(\gotn,0)\}$,
if $F_{\gotn,0}=[(X^{(1)}(\cdot))^{\gotn}]_{\vzero}\neq0$,
\item or the only segment $\{(0,\gotn),(\gotn-1,1)\}$,
if $F_{k,0}=0$ for all $k\ge0$,
\item or the two segments $\{(0,\gotn),(\gotn-1,1)\}$
and $\{(\gotn-1,1),(\gotp,0)\}$, with $\gotp \ge \gotn+1$,
if $F_{k,0}=0$ for all $k\le \gotp-1$ and $F_{\gotp,0}\neq0$.
\end{enumerate}

Therefore, in principle, a solution to the bifurcation equation
$F(\e,c)=0$ is
\begin{itemize}
\item{} either $c=c_{0}+\zeta_{1}\e +o(\e)$ for
some $\zeta_{1}\neq0$, corresponding to the segment
$\{(0,\gotn),(\gotn,0)\}$ if $[(X^{(1)}(\cdot))^{\gotn}]_{\vzero}\neq0$
(case 1) and to the segment $\{(0,\gotn),(\gotn-1,1)\}$ if
$[(X^{(1)}(\cdot))^{\gotn}]_{\vzero}=0$ (cases 2 and 3),
\item{} or $c=c_{0}+\zeta_{0}
\e^{\gotp-\gotn+1}+o(\e^{\gotp-\gotn+1})$ for some $\zeta_{0} \neq 0$,
corresponding to the segment $\{(\gotn-1,1),(\gotp,0)\}$ (case 3),
\item{} or $c=c_{0}+\zeta(\e)$, with $\zeta(\e)$ decaying to zero
faster than any power of $\e$ (case 2).
\end{itemize}

The case $\gotn=1$ should be discussed apart, but since it is
much easier and moreover has already been discussed in \cite{GBD1,GBD2},
here we concentrate on $\gotn>1$.

\begin{lemma} \label{lem:2.2}
If $[(X^{(1)}(\cdot))^{\gotn}]_{\vzero}=0$ then the bifurcation equation
(\ref{eq:2.5}) admits one and only one solution $c(\e)=o(\e)$.
Moreover such a solution is smooth in $\e$, and
it can be explicitly computed.
\end{lemma}

\prova If $[(X^{(1)}(\cdot))^{\gotn}]_{\vzero}=0$ one has
\begin{eqnarray}
\left[ (\zeta + \e X^{(1)}(\cdot) )^{\gotn} \right]_{\vzero}
& \!\!\! = \!\!\! &
\left[ (\zeta + \e X^{(1)}(\cdot) )^{\gotn} \right]_{\vzero} -
\left[ ( \e X^{(1)}(\cdot) )^{\gotn} \right]_{\vzero} =
\left[ (\zeta + \e X^{(1)}(\cdot) )^{\gotn} -
(\e X^{(1)}(\cdot))^{\gotn} \right]_{\vzero}
\nonumber \\
& \!\!\! = \!\!\! & \gotn \zeta \left[ \int_{0}^{1} {\rm d}s
\left( \e X^{(1)}(\cdot) + s\zeta \right)^{\gotn -1} \right]_{\vzero} =
\gotn \zeta \int_{0}^{1} {\rm d}s
\left[ \left( \e X^{(1)}(\cdot) + s\zeta \right)^{\gotn -1}
\right]_{\vzero} ,
\nonumber \end{eqnarray}
which is non-zero because $\gotn-1$ is even. Therefore,
by using that $\zeta=O(\e)$ by Lemma \ref{lem:2.1}, one obtains
\begin{equation}
F(\e,c) = \zeta \e^{\gotn-1} a +
G(\e,c) , \qquad a = g_{0} \gotn 
\int_{0}^{1} {\rm d}s
\left[ \left( X^{(1)}(\cdot) + s\zeta\e^{-1} \right)^{\gotn -1}
\right]_{\vzero} ,
\nonumber \end{equation}
with $a>0$ and $G(\e,c) =O(\e^{\gotn+1})$. Furthermore either
all derivatives of $G(\e,c_{0})$ vanish at $\e=0$
(and hence $F_{k,0}=0$ for all $k\ge 0$) or
\begin{equation}
G(\e,c_{0}) = b \e^{\gotp} + O(\e^{\gotp+1}) , \qquad
b = F_{\gotp,0} ,
\nonumber \end{equation}
for some $\gotp >\gotn$.

In both cases, there is no solution corresponding to the segment
$\{(0,\gotn),(\gotn-1,1)\}$. Indeed, the bifurcation
equation can be written as
\begin{equation}
\e^{\gotn-1} \left( a_{0} \zeta  + \Gamma(\e,\zeta) \right) = 0 ,
\qquad a_{0} = \left[ \left( X^{(1)}(\cdot) \right)^{\gotn -1}
\right]_{\vzero} ,
\nonumber \end{equation}
with the function $\Gamma(\e,\zeta)$ which is $C^{\io}$ in both $\e$
and $\zeta$. Hence we can apply the implicit function theorem to deduce
that for $\e\neq0$ there is one and only one solution $c=c_{0}+\zeta(\e)$
to the bifurcation equation, with $\zeta(\e)$ smooth in $\e$
and such that $\zeta(\e)=G(\e,c_{0})/a_{0}\e^{\gotn-1}+
\tilde\zeta(\e)$ and $\tilde\zeta(\e)/\zeta(\e)\to0$ as $\e\to0$.
Moreover one can explicitly estimate the interval $U$
around $\e=0$ in which such a solution exists (again by using
the implicit function theorem and the bounds in \cite{G2})).\EP

If $F_{k,0}=0$ for all $k\ge0$ then the only solution
to the bifurcation equation is $c=c_{0}+\zeta(\e)$,
where $\zeta(\e)$ is a $C^{\io}$ function which goes to zero
faster than any power as $\e$ tends to $0$.
If there exists $\gotp\ge\gotn+1$
such that $F_{k,0}=0$ for $k\le \gotp$ and $F_{p,0}\neq0$, the
only solution is of the form $c=c_{0} + \zeta_{0}\e^{\gotp-\gotn+1} +
o(\e^{\gotp-\gotn+1})$, where $\zeta_{0}$ is
the (unique) solution of the equation $a \zeta + b = 0$.

In the case $[(X^{(1)}(\cdot))^{\gotn}]_{\vzero}\neq0$
the following result holds.

\begin{lemma} \label{lem:2.3}
If $[(X^{(1)}(\cdot))^{\gotn}]_{\vzero}\neq 0$ then all solutions
of the bifurcation equation (\ref{eq:2.5}) are of order
$\e$ and have multiplicity 1. In particular, there exists a
constant $a_{0}>0$ such that $|c_{i}(\e)-c_{j}(\e)|>a_{0}|\e|$
for all pairs of such solutions $c_{i}(\e)$ and $c_{j}(\e)$.
\end{lemma}

\prova If $[(X^{(1)}(\cdot))^{\gotn}]_{\vzero}\neq0$ then the equation
\begin{equation}
\left[ \left( \zeta + \e X^{(1)}(\cdot)
\right)^{\gotn} \right]_{\vzero} = 0 ,
\label{eq:2.8} \end{equation}
with $X^{(1)}$ defined in (\ref{eq:2.4}), admits at least one
non-zero real solution $\zeta_{1}$ of order $\e$. Thus, one can write
\begin{equation}
\left( \zeta + \e X^{(1)}(\cdot) \right)^{\gotn} =
\left( \zeta - \zeta_{1} + \zeta_{1} + 
\e X^{(1)}(\cdot) \right)^{\gotn} =
\sum_{k=0}^{\gotn} \left( \begin{matrix} \gotn \\ k \end{matrix} \right)
\left( \zeta - \zeta_{1} \right)^{\gotn-k}
\left( \zeta_{1} + \e X^{(1)}(\cdot) \right)^{k} ,
\nonumber \end{equation} 
so that, by using that $[(\zeta_{1} + \e X^{(1)}
(\cdot))^{\gotn}]_{\vzero} = 0$, one has
\begin{equation}
\left[ \left( \zeta + \e X^{(1)}(\cdot) \right)^{\gotn}
\right]_{\vzero} =
\left( \zeta - \zeta_{1} \right)
\sum_{k=0}^{\gotn-1} 
\left( \begin{matrix} \gotn \\ k \end{matrix} \right)
\left( \zeta - \zeta_{1} \right)^{\gotn-1-k}
\left[ \left( \zeta_{1} + \e X^{(1)}(\cdot)
\right)^{k} \right]_{\vzero} .
\nonumber \end{equation} 
For $\zeta=\zeta_{1}$ one has
\begin{equation}
\sum_{k=0}^{\gotn-1}
\left( \begin{matrix} \gotn \\ k \end{matrix} \right)
\left( \zeta - \zeta_{1} \right)^{\gotn-1-k}
\left[ \left( \zeta_{1} + \e X^{(1)}(\cdot)
\right)^{k} \right]_{\vzero} =
\left( \begin{matrix} \gotn \\ \gotn-1 \end{matrix} \right)
\left[ \left( \zeta_{1} + \e X^{(1)}(\cdot) \right)^{\gotn-1}
\right]_{\vzero} > 0 ,
\nonumber \end{equation}
which shows that $\zeta_{1}$ is a simple root of the
equation (\ref{eq:2.8}), and hence also of the
bifurcation equation (\ref{eq:2.5}).\EP

The following result extends Lemma \ref{lem:2.2}, and
shows that the bifurcation equation admits a unique
solution for any $f$.

\begin{lemma} \label{lem:2.4}
The bifurcation equation (\ref{eq:2.5}) admits one and only one
solution $c(\e)$. One has $c(\e)=\zeta_{1}\e+O(\e^{2})$,
with $\zeta_{1}\neq0$ if $[(X^{(1)}(\cdot))^{\gotn}]_{\vzero}\neq0$
and $c(\e)=O(\e^{2})$ if $[(X^{(1)}(\cdot))^{\gotn}]_{\vzero}=0$.
\end{lemma}

\prova By Lemma \ref{lem:2.2} we can confine ourselves to the case
$[(X^{(1)}(\cdot))^{\gotn}]_{\vzero}\neq0$. In that case we know
by Lemma \ref{lem:2.3} that all the solutions of the bifurcation
equation are of order $\e$ and separated by order $\e$. Hence
we can set $\zeta=\zeta_{1}\e + \zeta'$, with $\zeta'=o(\e)$,
so that one has
\begin{eqnarray}
\frac{\partial}{\partial c} F(\e,c_{0}+\zeta)
& \!\!\! = \!\!\! &
g_{0} \frac{\partial}{\partial \zeta}
\left[ \left( \zeta + \e X^{(1)}(\cdot) \right)^{\gotn}
\right]_{\vzero} + O(\e^{\gotn+1}) \nonumber \\
& \!\!\! = \!\!\! &
g_{0}  \e^{\gotn} \gotn \left[ \left( \zeta_{1} + X^{(1)}(\cdot)
\right)^{\gotn-1} \right]_{\vzero} + O( \e^{\gotn+1}) ,
\nonumber \end{eqnarray}
which is strictly positive for $\e>0$ small enough. Thus,
at fixed $\e$, $F(\e,c)$ is increasing in $c$: this yields
that the solution $c=c(\e)$ of $F(\e,c)=0$ is unique for $\e$
small enough.\EP

\zerarcounters
\section{Multiscale analysis and diagrammatic expansion}
\label{sec:3}

In this section we develop a diagrammatic representation for the
quasi-periodic solution. There are a few differences with respect
to \cite{G2}, as we expand simultaneously both $X$ and $c$ in terms
of the perturbation parameter $\e$.
As in \cite{G2}, the expansions we shall find are not power
series expansions. See also \cite{GGG,CG} for analogous situations;
note, however, that in our case, as in \cite{GCB},
no fractional expansions arise.

We can summarise the results of the previous section as follows.
The bifurcation equation (\ref{eq:2.5}) admits a unique
solution $c(\e)=c_{0}+\zeta(\e)$ such that $\zeta(0)=0$.
Let $X^{(1)}$ the zero-average function such that $\dot X^{(1)}=f$.
If $[(X^{(1)}(\cdot))^{\gotn}]_{\vzero}\neq0$ then one has
$\zeta(\e)=\zeta_{1}\e + o(\e)$, where $\e\zeta_{1}$ is
the unique solution of the equation (\ref{eq:2.8}).

To make notation uniform we can set in the following $\zeta_{1}=0$
if $[(X^{(1)}(\cdot))^{\gotn}]_{\vzero}=0$. By defining
\begin{equation}
P_{\gotn}(\zeta)= \left[ \left( \zeta + X^{(1)}(\cdot) \right)^{\gotn}
\right]_{\vzero} , \qquad
P_{\gotn}'(\zeta) = \frac{{\rm d}}{{\rm d}\zeta} P_{\gotn}(\zeta) ,
\label{eq:3.0} \end{equation}
if $[(X^{(1)}(\cdot))^{\gotn}]_{\vzero}\neq0$
one has $P_{\gotn}'(\zeta_{1}) \neq 0$ by Lemma \ref{lem:2.3},
and if $[(X^{(1)}(\cdot))^{\gotn}]_{\vzero}=0$ one has
$P_{\gotn}'(\zeta_{1})=P_{\gotn}'(0) =\gotn
[(X^{(1)}(\cdot))^{\gotn-1}]_{\vzero}\neq0$ because $\gotn$ is odd.
Set $a=g_{0}P_{\gotn}'(\zeta_{1})$ in both cases; then $a \neq 0$.

In the following we give the diagrammatic rules in detail to
make the exposition self-contained, and we only stress where the
main differences lie with respect to the expansion of \cite{G2}.
From a technical point of view, the diagrammatic analysis
turns out a bit more involved, as it requires further expansions,
and hence more labels to be assigned to the diagrams.
On the other hand, eventually there is the advantage that one has a
completely constructive algorithm to determine the solution
within any fixed accuracy and an explicitly computable value
for the maximal size of the allowed perturbation.

A graph is a connected set of points and lines.
A \textit{tree} $\theta$ is a graph with no cycle,
such that all the lines are oriented toward a unique
point (\textit{root}) which has only one incident line (root line).
All the points in a tree except the root are called \textit{nodes}.
The orientation of the lines in a tree induces a partial ordering 
relation ($\preceq$) between the nodes. Given two nodes $v$ and $w$,
we shall write $w \prec v$ every time $v$ is along the path
(of lines) which connects $w$ to the root.

We call $E(\theta)$ the set of \textit{end nodes} in $\theta$,
that is the nodes which have no entering line, and $V(\theta)$
the set of \textit{internal nodes} in $\theta$, that is the set of
nodes which have at least one entering line. Set $N(\theta)=
E(\theta) \amalg V(\theta)$. With each end node $v$ we associate
a \textit{mode} label $\nn_{v}\in\ZZZ^{d}$. For all $v\in N(\theta)$
denote with $s_{v}$ the number of lines entering the node $v$;
for $v\in V(\theta)$ one has $s_{v} \ge \gotn$.

With respect to \cite{G2} the mode label of the end nodes can be
$\vzero$. This only occurs when $\zeta_{1}\neq0$. Hence we define
$E_{0}(\theta)=\{v\in E(\theta) : \nn_{v}=\vzero\}$ and
$E_{1}(\theta)=\{v\in E(\theta) : \nn_{v} \neq \vzero\}$;
we can set $E_{0}(\theta)=\emptyset$ if $\zeta_{1}=0$.
Define also $L_{0}(\theta)=\{\ell\in L(\theta) :
\ell \hbox{ exits a node } v\in E_{0}(\theta)\}$.
If $|N(\theta)|=1$ one requires $E_{1}(\theta)=N(\theta)$
and hence $E_{0}(\theta)=\emptyset$.

We denote with $L(\theta)$ the set of lines in $\theta$. Since a
line $\ell$ is uniquely identified with the node $v$ which it leaves,
we may write $\ell = \ell_{v}$. With each line $\ell$ we associate
a \textit{momentum} label $\nn_{\ell} \in \ZZZ^{d}$ and
a \textit{scale} label $n_{\ell}\in\ZZZ_{+} \cup \{-1\}$.
We set $n_{\ell}=-1$ when $\nn_{\ell}=\vzero$ (note that
$\nn_{\ell}=\vzero$ was not allowed in \cite{G2}).

The modes of the end nodes and the momenta of the lines
are related as follows: if $\ell = \ell_{v}$ one has
\begin{equation}
\nn_{\ell} = \sum_{w \in E(\theta) : w \preceq v} \nn_{w} .
\label{eq:3.1} \end{equation}
Given a tree $\theta$ we set
$\Lambda_{0}(\theta)=
\{\ell\in L(\theta): \nn_{\ell}=\vzero\} \setminus L_{0}(\theta)
= \{\ell\in L(\theta): \nn_{\ell}=\vzero$ and
$\ell$ does not exit an end node$\}$, and
define the \textit{order} of $\theta$ as
\begin{equation}
k(\theta) = |N(\theta)| - \gotn \left| \Lambda_{0}(\theta) \right| ,
\label{eq:3.2} \end{equation}
and the \textit{total momentum} of $\theta$ as $\nn(\theta)=
\nn_{\ell_{0}}$, if $\ell_{0}$ is the root line of $\theta$.

We call \textit{equivalent} two trees which can be transformed into
each other by continuously deforming the lines in such a way that
they do not cross each other. Let $\TT_{k,\nn}$ be the set of
inequivalent trees of order $k$ and total momentum $\nn$,
that is the set of inequivalent trees $\theta$ such that
$k(\theta)=k$ and $\nn(\theta)=\nn$.

A \textit{cluster} $T$ on scale $n$ is a maximal set of nodes and
lines connecting them such that all the lines have scales $n'\le n$
and there is at least one line with scale $n$. The lines entering
the cluster $T$ and the possible line coming out from it (unique if
existing at all) are called the external lines of the cluster $T$.
Given a cluster $T$ on scale $n$, we shall denote by $n_{T}=n$ the
scale of the cluster. We call $V(T)$, $E(T)$, and $L(T)$ the set of
internal nodes, of end nodes, and of lines of $T$, respectively;
the external lines of $T$ do not belong to $L(T)$.

We call \textit{self-energy cluster} any cluster $T$ such that
$T$ has only one entering line $\ell_{T}^{2}$ and one exiting
line $\ell_{T}^{1}$, and one has $\sum_{v\in E(T)} \nn_{v} = \vzero$
(and hence $\nn_{\ell_{T}^{1}}=\nn_{\ell_{T}^{2}}$ by (\ref{eq:3.1})).
Set $x_{T}=\oo\cdot\nn_{\ell_{T}^{2}}$.

We call \textit{excluded} a node $v$ such that $\nn_{\ell_{v}}=
\vzero$, $s_{v}=\gotn$, at least $\gotn-1$ lines entering $v$
do exit end nodes, and the other line $\ell'$ entering $v$
either also exits an end node or has momentum $\nn_{\ell'}=\vzero$.

Let $\gotT_{k,\nn}$ be the set of \textit{renormalised trees} in
$\TT_{k,\nn}$, i.e. of trees in $\TT_{k,\nn}$ which contain
neither any self-energy clusters nor any excluded nodes.
Define also $\gotR_{n}$ as the set
of renormalised self-energy clusters on scale $n$, i.e.
of self-energy clusters on scale $n$ which contain neither any
further self-energy clusters nor any excluded nodes.

\begin{lemma} \label{lem:3.1}
One has $k(\theta) \ge 1$ if $\nn(\theta)\neq\vzero$
and $k(\theta) \ge 2$ if $\nn(\theta)=\vzero$.
\end{lemma}

\prova By induction on $|N(\theta)|$. If $|N(\theta)|=1$ then
$N(\theta)=E_{1}(\theta)$ and hence $\Lambda_{0}(\theta)=\emptyset$,
so that $k(\theta)=1$. If $|N(\theta)|>1$ let $v_{0}$ be the last node
of $\theta$, that is the node which the root line of $\theta$ exits,
and let $\theta_{1},\ldots,\theta_{s}$ the trees with root in $v_{0}$.
Note that $s\ge \gotn$. Then one has
$|N(\theta)|=1+|N(\theta_{1})|+\ldots+|N(\theta_{s})|$, while
$|\Lambda_{0}(\theta)|=|\Lambda_{0}(\theta_{1})|+\ldots
+|\Lambda_{0}(\theta_{s})|$ if $\nn(\theta)\neq\vzero$ and
$|\Lambda_{0}(\theta)|=1+|\Lambda_{0}(\theta_{1})|+\ldots
+|\Lambda_{0}(\theta_{s})|$ if $\nn(\theta)=\vzero$.
In the first case one has $k=k(\theta)=1+k(\theta_{1})+\ldots+
k(\theta_{s})\ge 1+\gotn$. In the second case one has
$k=k(\theta)=1+k(\theta_{1})+\ldots+
k(\theta_{s}) - \gotn \ge 1+s-\gotn$, so that $k\ge1$.
Moreover $k=1$ would be possible only if $s=\gotn$
and $k(\theta_{1})=\ldots=k(\theta_{s})=1$. However, in such a case
the lines entering $v_{0}$ would all exit end nodes and hence
$v_{0}$ would be an excluded node. Thus, $k\ge2$ if
$\nn(\theta)=\vzero$.\EP

\begin{lemma} \label{lem:3.2}
There exists a positive constant $\kappa$ such that $|N(\theta)| \le
\kappa k(\theta)$ for any renormalised tree $\theta$.
One can take $\kappa=3\gotn$.
\end{lemma}

\prova We give the proof for $\gotn>1$
(the case $\gotn=1$ being much easier; see \cite{GBD1}).

For $k(\theta)=1$ the bound is trivially satisfied.
We prove that for all $k\ge 2$, all $\nn\in\ZZZ^{d}$ and
all trees $\theta\in\gotT_{k,\nn}$ one has  $|N(\theta)| \le
a k(\theta)-b$, with $a=3\gotn$ and $b=4\gotn$.

The proof is by induction on $k$. For $k=2$ it is just a check:
if $\nn\ne\vzero$ one has $N(\theta)=2$, while if $\nn=\vzero$
one has $N(\theta)=2+\gotn\le 2\gotn$ for $\gotn>1$.

By assuming that the bound holds for all $k< h$, one has
for $\theta\in\gotT_{h,\nn}$
\begin{equation}
|N(\theta)| = s_{1} + \sum_{i=1}^{s_{0}} |N(\theta_{i})| ,
\qquad |\Lambda_{0}(\theta)| \le 1 +
\sum_{i=1}^{s_{0}} |\Lambda_{0}(\theta_{i})| ,
\nonumber \end{equation}
where $\ell_{1},\ldots,\ell_{s_{0}}$ are the lines
in $\Lambda_{0}(\theta)$ closest to the root line of $\theta$,
$\theta_{1},\ldots,\theta_{s_{0}}$ are the trees with root lines
$\ell_{1},\ldots,\ell_{s_{0}}$, respectively, and $s_{1}$ is the
number of nodes which precede the root line of $\theta$ but not
the root lines of $\theta_{1},\ldots,\theta_{s_{0}}$.
By using the definition of order (\ref{eq:3.2}) one has
\begin{equation}
h = k(\theta) \ge s_{1} - \gotn + \sum_{i=1}^{s_{0}}
k(\theta_{i}) .
\nonumber \end{equation}
Note that $k(\theta_{i})\ge 2$ by Lemma \ref{lem:3.2}, and hence
$|N(\theta_{i})| \le a k(\theta_{i})-b$ by the inductive hypothesis.

If $s_{0}=0$ then $\Lambda_{0}(\theta) \le 1$, so that
$|N(\theta)| \le k(\theta) + \gotn = h + \gotn
\le 3\gotn h - 4\gotn$ for $\gotn>1$ and $h>1$.

If $s_{0} \ge 1$ the inductive hypothesis yields
\begin{eqnarray}
|N(\theta)| & \!\!\! \le \!\!\! &
s_{1} + a \sum_{i=1}^{s_{0}} k(\theta_{i}) -
s_{0} b \nonumber \\
& \!\!\! \le \!\!\! &
a k(\theta) - b - 
\Big( (a-1) s_{1} + ( s_{0} - 1 ) b - a\gotn \Big) .
\nonumber \end{eqnarray}
Thus, the assertion follows if
\begin{equation}
\left( a-1 \right) s_{1} + \left( s_{0} - 1 \right) b -
a \gotn = a \left( s_{1} + s_{0} - \gotn \right) +
\left( b-a \right) \left( s_{0}-1 \right) - s_{1} - a  \ge 0 .
\nonumber \end{equation}

A key remark is that $s_{0}+s_{1} = \gotn + 1 + p$,
with $p \ge 0$, because the root line of $\theta$ exits
a node $v_{0}\in V(\theta)$, so that $s_{v_{0}} \ge \gotn$.
Then we can rewrite
\begin{equation}
a \left( s_{1} + s_{0} - \gotn \right) +
\left( b-a \right) \left( s_{0}-1 \right) - s_{1} - a =
a p - \gotn - p - 1 +
s_{0} + \left( b - a \right) \left( s_{0}-1 \right) .
\nonumber \end{equation}
If $p=0$ then $s_{1}+s_{0}=\gotn +1$; this yields $s_{0} \ge 2$
(otherwise $v_{0}$ would be an excluded node), so that
$a p - \gotn - p - 1 +
s_{0} + \left( b - a \right) \left( s_{0}-1 \right) \ge
1 + \left( b - a \right) - \gotn \ge 1$.
If $p \ge 1$ then $a p - \gotn - p - 1 +
s_{0} + \left( b - a \right) \left( s_{0}-1 \right) \ge
a p - \gotn - p = \gotn p + \gotn \left( p - 1 \right) +
p \left( \gotn - 1 \right) \ge 2\gotn-1$.\EP


Let $\psi$ be a non-decreasing $C^{\infty}$ function defined
in $\RRR_{+}$, such that
\begin{equation}
\psi(u) = \left\{
\begin{array}{ll}
1 , & \text{for } u \geq 1 , \\
0 , & \text{for } u \leq 1/2 ,
\end{array} \right.
\label{eq:3.3} \end{equation}
and set $\chi(u) := 1-\psi(u)$. For all $n \in \ZZZ_{+}$ define
$\chi_{n}(u) := \chi(u/4\al_{n}(\oo))$ and $\psi_{n}(u) :=
\psi(u/4\al_{n}(\oo))$, and set
$\Xi_{0}(x)=\chi_{0}(|x|)$, $\Psi_{0}(x)=\psi_{0}(|x|)$, and
\begin{equation}
\Xi_{n}(x)=\chi_{0}(|x|)\ldots \chi_{n-1}(|x|) \chi_{n}(|x|) , \qquad
\Psi_{n}(x)=\chi_{0}(|x|)\ldots \chi_{n-1}(|x|) \psi_{n}(|x|) ,
\label{eq:3.4} \end{equation}
for $n\ge 1$. Then we define the node factor as
\begin{equation}
F_{v} = \begin{cases}
- \displaystyle{ \frac{1}{s_{v}!}
\frac{{\rm d}^{s_{v}}}{{\rm d}x^{s_{v}}} g(c_{0}) } ,
& v \in V(\theta) , \\
f_{\nn_{v}} , & v \in E_{1}(\theta) , \\
\zeta_{1} , & v \in E_{0}(\theta) ,
\end{cases}
\label{eq:3.5} \end{equation}
and the propagator as
\begin{equation}
G_{\ell} = \begin{cases}
G^{[n_{\ell}]}(\oo\cdot\nn_{\ell};\e,c_{0}) , & \nn_{\ell} \neq \vzero , \\
- \displaystyle{ \frac{1}{a} } , & \nn_{\ell} = \vzero ,
\quad \ell \in \Lambda_{0}(\theta) , \\
1 , & \nn_{\ell} = \vzero , \quad \ell \in L_{0}(\theta) , 
\end{cases}
\label{eq:3.6} \end{equation}
with $\zeta_{1}$ and $a$ defined before and after (\ref{eq:3.0}),
respectively, and $G^{[n]}(x;\e,c_{0})$ recursively
defined for $n\ge 0$ as
\begin{subequations}
\begin{align}
& \hskip-.3truecm
G^{[n]}(x;\e,c) =
\frac{\Psi_{n}(x)}{ \ii x(1+\ii \e x) - \MM^{[n-1]}(x;\e,c)} ,
\label{eq:3.7a} \\
& \hskip-.3truecm
\MM^{[n]}(x;\e,c) = \MM^{[n-1]}(x;\e,c) +
\Xi_{n}(x)
M^{[n]}(x;\e,c) , \quad
M^{[n]}(x;\e,c) = \!\! \sum_{T\in\gotR_{n}}\Val(T,x;\e,c) ,
\label{eq:3.7b}
\end{align}
\label{eq:3.7} \end{subequations}
where $\MM^{[-1]}(x;\e,c_{0})=\e g_{1}(c_{0})=0$ and
\begin{equation}
\Val(T,x;\e,c_{0}) = \Big(\prod_{\ell \in L(T)} G_{\ell} \Big) 
\Big( \prod_{v \in N(T)} F_{v} \Big)
\label{eq:3.8} \end{equation}
is called the value of the self-energy cluster $T$.

Note that, with respect to the analogous formulae (3.5) and (3.6)
of \cite{G2}, both $F_{v}$ and $G_{\ell}$ are computed at $c=c_{0}$.
We prefer writing explicitly the dependence on $c_{0}$
(even if not necessary, as $c_{0}$ is fixed once and for all
as in Assumption \ref{hyp:2}) simply in order to use
the same notations as in \cite{G2}.

Set
\begin{equation}
X^{[k]}_{\nn} = \sum_{\theta\in \gotT_{k,\nn}}
\Val(\theta;\e,c_{0}) ,
\qquad \nn \neq \vzero, \quad k \ge 1 ,
\label{eq:3.9} \end{equation}
and
\begin{equation}
\zeta^{[k]}= \sum_{\theta\in \gotT_{k+\gotn,\vzero}}
\Val(\theta;\e,c_{0}) ,
\qquad k \ge 2 ,
\label{eq:3.10} \end{equation}
where the tree value $\Val(\theta;\e,c_{0})$ is defined as
\begin{equation}
\Val(\theta,x;\e,c_{0}) = \Big(\prod_{\ell \in L(\theta)} G_{\ell} \Big) 
\Big( \prod_{v \in N(\theta)} F_{v} \Big) .
\label{eq:3.11} \end{equation}

Call $\gotN_{n}(\theta)$ the number of lines $\ell\in L(\theta)$
such that $n_{\ell}\ge n$, and $\gotN_{n}(T)$ the number of lines
$\ell\in L(T)$ such that $n_{\ell}\ge n$, and set
\begin{equation}
M(\theta) = \sum_{v \in E(\theta)} |\nn_{v}| , \qquad
M(T) = \sum_{v \in E(T)} |\nn_{v}| .
\label{eq:3.12} \end{equation}
%
%
%

Define the renormalised series
\begin{equation}
\overline X(\pps;\e)= \sum_{\nn\in\ZZZ^{d}_{*}}
{\rm e}^{\ii\nn\cdot\pps} \overline X_{\nn} ,
\qquad \overline X_{\nn} = \sum_{k=1}^{\io} \e^{k} X_{\nn}^{[k]} .
\label{eq:3.14} \end{equation}
and
\begin{equation}
\overline \zeta(\e) = \e \zeta_{1} + \sum_{k=2}^{\io} \e^{k} \zeta^{[k]} .
\label{eq:3.15} \end{equation}
In the next section we shall prove first that
the series (\ref{eq:3.14}) and (\ref{eq:3.15}) converge,
then that the function $x_{0}(t) = c_{0} + \overline \zeta(\e) +
\overline X(\oo t;\e)$ solves the equations
\begin{subequations}
\begin{align}
& \ii\oo\cdot\nn \left( 1 + \ii\e\oo\cdot\nn \right) X_{\nn} +
\e \left[ g(c_{0}+\zeta(\e)+X(\cdot;\e) \right]_{\nn} = \e f_{\nn} ,
\qquad \nn \neq \vzero ,
\label{eq:3.16a} \\
& [g(c_{0}+\zeta(\e)+X(\cdot;\e)]_{\vzero} - f_{\vzero} = 0 ,
\label{eq:3.16b}
\end{align}
\label{eq:3.16} \end{subequations}
\vskip-.3truecm
\noindent and hence the equation (\ref{eq:1.1}).

\zerarcounters
\section{Bounds}
\label{sec:4}

The proof of the convergence of the renormalised series
proceeds as in \cite{G2}. We confine ourselves to state the
basic steps of the proof, without giving the details,
except when the discussion departs from \cite{G2}.

\begin{lemma} \label{lem:4.1}
For any renormalised tree $\theta$,
one has $\gotN_{n}(\theta) \le 2^{-(n-2)}M(\theta)$.
\end{lemma}

\prova The same as the proof of Lemma 3.1 in \cite{G2}.\EP

\begin{lemma} \label{lem:4.2}
Assume there exists a constant $C_{0}$ such that
$|G^{[n]}(x;\e,c_{0})|\le C_{0}/\al_{n}(\oo)$ for all $n\in\ZZZ_{+}$.
Then there exists $\e_{0}>0$ such that,
for all all $|\e|<\e_{0}$,
the series (\ref{eq:3.14}) and (\ref{eq:3.15}) converge.
Moreover the series (\ref{eq:3.14}) is analytic in $\pps$.
\end{lemma}

\prova Let $\theta$ be a tree in $\gotT_{k,\nn}$ and set
$K=|N(\theta)|$. One reasons as in \cite{G2} to prove that
\begin{equation}
\left| \Val(\theta;\e,c) \right| \le 
\widetilde C_{0}^{K} D_{0}^{K} \al_{n_{0}}^{-K}(\oo)
{\rm e}^{-\xi M(\theta)/2} ,
\nonumber \end{equation}
where $D_{0}=\max\{\Gamma,\Phi\}$, $\widetilde C_{0}=\max\{C_{0},a\}$,
and $n_{0} \in \NNN$ such that
\begin{equation}
4 \sum_{n=n_{0}+1}^{\io} \frac{1}{2^{n}}
\log \frac{1}{\al_{n}(\oo)} \le \frac{\xi}{2} .
\nonumber \end{equation}
By Lemma \ref{lem:3.2} one has $K\le \kappa k$.
The number of trees of order $k$ with fixed mode labels is hence
bounded by $C_{2}^{k}$ for a suitable constant $C_{2}$, and the sum
over the mode labels can be performed by using half the exponent
in the dacaying factor ${\rm e}^{-\xi M(\theta)/2}$. This gives
\begin{equation}
\left| X^{[k]}_{\nn} \right| \le C^{k} {\rm e}^{-\xi|\nn|/4} ,
\nonumber \end{equation}
for some constant $C$,
which yields the convergence of the series (\ref{eq:3.14})
to a function analytic in $\pps$. In the same way one obtains
\begin{equation}
\left| \zeta^{[k]} \right| \le C^{k} 
\nonumber \end{equation}
with the same constant $C$, and
this yields the convergence of the series (\ref{eq:3.15}).\EP

\begin{lemma} \label{lem:4.3}
For any self-energy cluster $T\in\gotR_{n}$ such that
$\Xi_{n}(x_{T}) \neq 0$, one has
$M(T) \ge 2^{n-1}$ and $\gotN_{p}(T) \le 2^{-(p-2)}M(T)$
for all $p \le n$.
\end{lemma}

\prova The same as the proof of Lemma 3.3 of \cite{G2}.\EP

\begin{lemma} \label{lem:4.4}
Assume the propagators $G^{[p]}(x;\e,c_{0})$ are differentiable in $x$
and there exist constants $C_{0}$ and $C_{1}$ such that
$|G^{[p]}(x;\e,c_{0})| \le C_{0}/\al_{p}(\oo)$ and
$|\partial_{x} G^{[p]}(x;\e,c_{0})| \le C_{1}/\al_{p}^{3}(\oo)$ for all $p<n$.
Then there exists $\e_{0}>0$ such that, for all $|\e|<\e_{0}$,
the function $x \mapsto M^{[n]}(x;\e,c_{0})$ is differentiable, and one has
\begin{equation}
\left| M^{[n]} (x;\e,c_{0}) \right| ,
\left| \partial_{x} M^{[n]} (x;\e,c_{0}) \right| \le
D_{1} |\e|^{2} {\rm e}^{-D_{2}2^{n}} ,
\nonumber \end{equation}
for some positive constants $D_{1}$ and $D_{2}$.
\end{lemma}

\prova Proceed as in the proof of Lemma 3.4 of \cite{G2},
by taking into account the differences already pointed out
in the proof of Lemma \ref{lem:4.2}.\EP

\begin{lemma} \label{lem:4.5}
Assume there exists a constant $C_{0}$ such that
$|G^{[p]}(x;\e,c_{0})|\le C_{0}/\al_{p}(\oo)$ for all $p<n$.
Then one has $(\MM^{[p]}(x;\e,c_{0}))^{*} = \MM^{[p]}(-x;\e,c_{0})$
for all $p \le n$.
\end{lemma}

\prova As the proof of Lemma 3.5 of \cite{G2}.\EP

\begin{lemma} \label{lem:4.6}
Let $\e_{0}$ as in Lemma \ref{lem:4.1}.
For all $n\in\ZZZ_{+}$ the function $x \mapsto \MM^{[n]}(x;\e,c_{0}))$
is differentiable and one has $|\ii x (1 + \ii \e x) -
\MM^{[n]}(x;\e,c_{0})|\ge |x|/2$ for all $|\e|<\e_{0}$.
\end{lemma}

\prova As the proof of Lemma 3.6 of \cite{G2}.\EP 

By collecting together the results above we can prove that
the series (\ref{eq:3.14}) and (\ref{eq:3.15}) converge.
The proof proceeds as follows. Convergence holds if the
hypotheses of Lemma \ref{lem:4.2} are verified. So we have to
prove that the propagators satisfy the bounds 
$|G^{[n]}(x;\e,c_{0})| \le C_{0}/\al_{n}(\om)$.
More precisely, we can prove by induction on $n$
that the propagators satisfy the bounds 
\begin{equation}
\left| G^{[n]}(x;\e,c_{0}) \right| \le \frac{C_{0}}{\al_{n}(\oo)} ,
\qquad \left| \partial_{x} G^{[n]}(x;\e,c_{0}) \right|
\le \frac{C_{1}}{\al_{n}^{3}(\oo)} ,
\label{eq:4.1} \end{equation}
for suitable constants $C_{0}$ and $C_{1}$. For $n=0$
the check is trivial. Then we assume that the bounds are satisfied
up to scale $n-1$, and we write $G^{[n]}(x;\e,c_{0})$ according
to (\ref{eq:3.7}). Now, $G^{[n]}(x;\e,c_{0})$ depends on the
quantities $\MM^{[p]}(x;\e,c_{0})$ with $p < n$,
which in turn depend on the propagators $G^{[p']}(x;\e,c_{0})$
with $p'\le p$. Thus, by the inductive
hypothesis we can apply Lemma \ref{lem:4.4} and Lemma \ref{lem:4.5}
to deduce Lemma \ref{lem:4.6}. This implies that also
$G^{[n]}(x;\e,c_{0})$ satisfies the bounds (\ref{eq:4.1}).
Smoothness in $\e$ of the series is discussed in the same way,
proving by induction that the derivatives of the propagators
satisfy the bounds $|\partial_{\e}^{m}G^{[n]}(x;\e,c_{0})| \le
K_{m}/\al_{n}^{m}(\oo)$ for all $m\ge1$ and for suitable
constants $K_{m}$ (we refer to \cite{G2} for further details).

Finally, by proceeding as in the proof of Lemma 3.7 of \cite{G2},
we can prove that the function $x_{0}(t)=c_{0}+\zeta(\e)+
X(\oo t;\e)$ solves the equations (\ref{eq:3.16}).
As far as equation (\ref{eq:3.16a}) is concerned the proof proceeds
as in the proof of Lemma 3.7 of \cite{G2}. To deal with
equation (\ref{eq:3.16b}), we write
$[g(c_{0}+\zeta(\e)+X(\cdot;\e)]_{\vzero}$ according to
(\ref{eq:2.6}) and expand $X(\cdot;\e)=\overline X(\cdot;\e)$
according to (\ref{eq:3.14}) and (\ref{eq:3.9}). Then we obtain
\begin{equation}
g_{0} P_{\gotn}'(\zeta_{1}) \left( \zeta_{k} + 
\sum_{\theta\in \gotT_{k+\gotn,\vzero}}
\Val(\theta;\e,c_{0}) , \right) \qquad k \ge 2 ,
\label{eq:4.2} \end{equation}
which is identically zero if $\zeta_{k}$ is defined
by (\ref{eq:3.10}).



\begin{thebibliography}{99}

{\small \bibitem{BDG}{
M.V. Bartuccelli, J.H.B. Deane, G. Gentile,
\textit{Globally and locally attractive solutions
for quasi-periodically forced systems},
J. Math. Anal. Appl.
\textbf{328} (2007), no. 1, 699-714. }
%
\bibitem{BK}{
E. Brieskorn, H. Kn\"orrer,
\textit{Plane algebraic curves},
Birkh\"auser, Basel, 1986. }
%
\bibitem{CH}{
L.O. Chua, A. Ushida,
\textit{Algorithms for computing almost periodic steady-state
response of nonlinear systems to multiple input frequencies},
IEEE Trans. Circuits and Systems
\textbf{28} (1981), no. 10, 953--971. } 
%
\bibitem{CG}{
L. Corsi, G. Gentile,
\textit{Melnikov theory to all orders and Puiseux series
for subharmonic solutions},
J. Math. Phys.
\textbf{49} (2008), no. 11, 112701, 29 pp. }
%
\bibitem{GGG}{
G. Gallavotti, G. Gentile, A. Giuliani,
\textit{Fractional Lindstedt series},
J. Math. Phys. \textbf{47} (2006), no. 1, 012702, 33 pp. }
%
\bibitem{G1}{
G. Gentile,
\textit{Quasi-periodic solutions for two-level systems},
Comm. Math. Phys.
\textbf{242} (2003),  no. 1-2, 221--250. }
%
\bibitem{GCB}{
G. Gentile, D.A. Cortez, J.C.A. Barata,
\textit{Stability for quasi-periodically perturbed Hill's equations},
Comm. Math. Phys.
\textbf{260} (2005), no. 2, 403--443. }
%
\bibitem{G2}{
G. Gentile,
\textit{Quasi-periodic motions in strongly dissipative forced systems},
Ergodic Theory Dynam. Systems, to appear. }
%
\bibitem{GBD1}{
G. Gentile, M.V. Bartuccelli, J.H.B. Deane,
\textit{Summation of divergent series and Borel summability
for strongly dissipative differential equations with periodic
or quasiperiodic forcing terms},
J. Math. Phys. \textbf{46} (2005), no. 6, 062704, 21 pp. }
%
\bibitem{GBD2}{
G. Gentile, M.V. Bartuccelli, J.H.B. Deane,
\textit{Quasiperiodic attractors, Borel summability
and the Bryuno condition for strongly dissipative systems},
J. Math. Phys. \textbf{47} (2006), no. 7, 072702, 10 pp. }
%
\bibitem{KP1}{
Ch. Kaas-Petersen,
\textit{Computation of quasiperiodic solutions
of forced dissipative systems},
J. Comput. Phys.
\textit{58} (1985), no. 3, 395--408. }
%
\bibitem{KP2}{
Ch. Kaas-Petersen,
\textit{Computation of quasiperiodic solutions
of forced dissipative systems. II},
J. Comput. Phys.
\textit{64} (1986), no. 3, 433--442. }

}

\end{thebibliography}
\end{document}